\documentstyle[amscd,12pt,psamsfonts,leqno,amssymb]{amsart}

\input xypic

\newtheorem{defn}[equation]{Definition}
\newtheorem{thm}[equation]{Theorem}

\newtheorem{cor}[equation]{Corollary}
\newtheorem{lem}[equation]{Lemma}
\newtheorem{prop}[equation]{Proposition}

\theoremstyle{remark}
\newtheorem{rem}{Remark}
\theoremstyle{remark}

\numberwithin{equation}{section}

\renewcommand\sp{\operatorname{Spec}}

\newcommand\grv{{\operatorname{Gr}}(V)}
\newcommand\gr{\operatorname{Gr}}

\renewcommand\hom{\operatorname{Hom}}
\renewcommand\det{\operatorname{det}}
\newcommand\Det{\operatorname{Det}}

\newcommand\limpl[1]{\underset{#1}\varprojlim\,}

\newcommand\aut{\operatorname{Aut}}

\newcommand\rk{\operatorname{rk}}

\renewcommand\o{{\mathcal O}}

\renewcommand\L{{\mathcal L}}
\renewcommand\c{{\mathcal C}}

\newcommand\Z{{\mathbb Z}}

\begin{document}

\title[Arithmetic Grassmannians and Central Extensions]{Arithmetic infinite Grassmannians and \\
the induced central extensions}
\author[F. J. Plaza Mart\'{\i}n]{Francisco J. Plaza Mart\'{\i}n}
   \address{Departamento de Matem\'aticas, Universidad de Salamanca,  Plaza
   de la Merced 1-4
   \\
   37008 Salamanca. Spain.
   \\
    Tel: +34 923294460. Fax: +34 923294583}

\thanks{
\indent This work is partially supported by the research contracts MTM2006-0768 of DGI
and SA112A07 of JCyL}
\email{fplaza@@usal.es}

\begin{abstract}
The construction of families of Sato Grassmannians, their determinant line bundles and the extensions induced by them are given. The base scheme is an arbitrary scheme.
\end{abstract}


\maketitle


\section{Introduction}

Since infinite Grassmannians were introduced by Sato (\cite{SS}) as a classifying space for the solutions of the KP hierarchy, their structures and properties have been successfully applied in a wide range of topics, such as microlocal analysis (\cite{DJKM}), loop groups (\cite{SW}), conformal and quantum field theories and string theory (\cite{W, KNTY, MP2}), moduli problems on algebraic curves (\cite{M, Shiota, MP}), representation theory of infinite dimensional lie algebras (\cite{Kac-Raina}), abelian and non-abelian reciprocity laws on curves (\cite{AP,MPa}), supersymmetric analogues (\cite{B,MR}), etc.

Some of the above-cited works are strongly based on the algebraic structure of Sato Grassmannians (\cite{AMP}). For instance, once a Sato Grassmannian has been endowed with an algebraic structure it makes sense to apply algebraic geometry to study line bundles on it and automorphisms of it or of certain line bundles. It is worth to mention that the functor of infinite dimensional Grassmannian has been studied over Grassmann algebras (\cite{B,MR}) as well as on the category of noetherian schemes over a field, see \cite{Quandt} (but, unfortunately, it fails to be representable in this category). However, up to now only the representability of Sato Grassmannians of suitable infinite dimensional $k$-vector spaces on the category of $k$-schemes, where $k$ is a field, were known.

This paper aims at offering a construction of Sato Grassmannians of certain $\o_S$-modules in the category of $S$-schemes where $S$ is an arbitrary scheme; i.e. the base field is no longer required. Note that, from the algebro-geometric point of view, this new object is the natural framework within which to study the variation of coefficients in the KP-hierarchy. Besides its interest in its own right, let us say few words on some issues that could profit from this construction. First, let us recall that some main characters of the geometric Langlands program, such as the moduli space of pointed curves, the moduli space of $G$-bundles, etc. have been constructed with the help of Sato Grassmannians over the field of complex numbers. Now, such constructions can be generalized for other base schemes (allowing families of moduli spaces, function fields as base schemes, etc.). Indeed, this new object has already been  applied in \cite{HMP-Higgspairs}, where $S$ parametrizes $k$-algebra structures in the $k$-vector space. For the interested reader, let us recommend  Frenkel's survey \cite{F}. Second, the Sato Grassmannian presented here would also provide a geometric formulation of the approach to reciprocity laws given in \cite{BBE} (see also \cite{AP} where the base is the spectrum of a local artinian ring). Note that the case of valuation rings of unequal characteristic is also allowed.

A second main result are the constructions of the determinant line bundle (recall that its space of global sections is deeply connected to Fock spaces \cite{Kac-Raina}) and of a central extension that includes the transformation groups considered in several works (e.g. \cite{AP,MPa,MP2},...) that will be suitable for explicit computations. Let us remark that our determinants and extensions are deeply related with those of \S2 of {\cite{BBE} and this relationship deserves further study. In the future, we shall apply these results to some moduli and arithmetic problems for families of arithmetic curves; more precisely, in the problem of a reciprocity law for families of arithmetic curves defined over ${\mathbb Z}$ and over $p$-adic numbers.

Let us briefly explain  how the paper is organized. Section~2 is devoted to the construction of the Sato Grassmannian of an $S$-module $V$ (under certain hypotheses) as an $S$-scheme. The main steps of this construction are inspired by \cite{AMP} and by Grothendieck's constructions of Grassmannians (\cite{EGAIHES}), although some arguments are now quite technical and delicate. This section ends with explicit computations of two spaces of global sections, which will be required later.

Section 3 focuses on the construction of extensions with the help of Sato Grassmannians. Contrary to the case of a base field, now we must deal with monoids rather than with groups. Thus, Leech's (\cite{Leech}) results on extension of groups by monoids will be extensively used. Indeed, once the determinant line bundle has been constructed (\S3.A), the central extension shows up as the set of certain transformations of it endowed with the composition law (\S3.B) and a criterion for this extension to be trivial is proven. Nevertheless, it turns out that, in the case of families (e.g. $S$ arbitrary), the defining conditions may be too strong and, consequently, the monoid too small. In \S3.C we propose a method to overcome this difficulty with the help of $K$-theory. Due to the close relationship between these two monoids, we expect that explicit computations can be carried out.

Finally, I wish to express my gratitude to Professors G Anderson, M Mulase and J M Mu\~{n}oz Porras for their valuable suggestions and comments.

\section{Arithmetic Sato Grassmannians}
\subsection{Commensurableness}

Henceforth, we consider an arbitrary scheme, $S$, and a pair of sheaves, $(V, V^+)$, of flat
quasicoherent $\o_S$-modules such that $V^+ \subset V$ and $V/V^+$ are flat.

\begin{defn}\label{defn:commen}
Let $A,B\subseteq V$ be quasicoherent $\o_S$-submodules such that
$V/A,V/B$ are flat.

We say that $A$ and $B$ are commensurable, $A\sim B$, if  $V/A+B$
and $V/A\cap B$ are flat and $A+B/A\cap B$ is locally free of finite
type (l.f.f.t.).
\end{defn}

Observe that the family $\{A\subseteq V \,\vert\, A\sim V^+\}$
canonically induces a topology on $V$ and also topologies in submodules
and quotients of $V$; these topologies will be called
$V^+$-topologies.

Furthermore, it will be assumed that the following three conditions hold
    \begin{itemize}
    \item $V=\bigcup_{A\sim V^+}A$;
    \item separateness $(0)=\bigcap_{A\sim V^+}A$;
    \item completeness $V=\limpl{A\sim V^+}V/A$.
    \end{itemize}
Observe that $V/A = \hat V/\hat A=\limpl{B\sim A} V/(A+B)$.

Let $T\to S$ be an $S$-scheme and $W$ a module endowed with the
$V^+$-topology. Then we write
    $$
    \hat W_T\,:=\, W \hat \otimes_{\o_S}\o_T
    $$
where $\hat\otimes$ denotes the completion w.r.t. the
$V^+$-topology. Note that for a morphism $R\to T$ of $S$-schemes one
has that $(\hat W_T)^{\hat{\,}}_R=\hat W_R$.

Given $(V,V^+)$ and an $S$-scheme, $T$,  we may construct a new pair
satisfying the above requirements; namely, $(\hat V_T, \hat V^+_T)$.

\begin{lem}\label{lem:lfftquotientimpliscomm}
Let $B\subseteq A\subseteq V$ be $\o_S$-submodules such that $A$ is
quasicoherent and $V/A$ is flat.

If $A/B$ is l.f.f.t., then $B$ is quasicoherent, $V/B$ is flat and
$A\sim B$.
\end{lem}

\begin{pf}
From the hypothesis, we know that $A/B$ and $V/A$ are quasicoherent
and flat. Therefore, the exactness of the sequence
    $$
    0\to A/B \to V/B \to V/A \to 0
    $$
implies that $V/B$ is quasicoherent and flat. Finally, $B$ is also
quasicoherent since it is the kernel of $V\to V/B$, which is a
morphism of quasicoherent sheaves.
\end{pf}

\begin{thm}
Commensurableness is an equivalence relation and it is preserved
under base change.
\end{thm}

\begin{pf}
to say that it is preserved under base change means that $A\sim B$ implies $\hat A_T\sim\hat  B_T$ for every morphism $T\to S$ and every $A,B$, as in definition~\ref{defn:commen}. This  is straightforward.

The first part of the statement will follow easily if we can prove that
    $$
    A\sim B \qquad \implies \qquad A\sim A\cap B\text{ and } A\sim A+ B
    $$

From the hypothesis we know that $A+B$ is quasicoherent and that
$V/A+B$ is flat. By Lemma~\ref{lem:lfftquotientimpliscomm} is
suffices to show that $A+B/A$ is l.f.f.t.~. Observe that the exact sequence
    $$
    0 \to A/A\cap B \to A+B/A\cap B \to A+B/A \to 0
    $$
shows that $A+B/A$ is  of finite type. Recall that $V, V/A, V/B,
V/A+B, V/A\cap B$ are flat and that, therefore, $A,B, A+B$ and
$A\cap B$ are flat, and hence so are the three terms in the
sequence. Using $A+B/A\simeq B/A\cap B\subseteq A+B/A\cap B$,
flatness and Nakayama's Lemma one has that $A+B/A$ is l.f.f.t.~.

Similarly, one proves that $A/A\cap B\simeq A+B/B$ is locally free
of finite type.
\end{pf}

\begin{rem}
Let $A,B$ be as in definition~\ref{defn:commen} and $A\sim B$.  Note that the subset of $S$ given by $\{s\in S\,\vert\, A_{k(s)}\subseteq B_{k(s)}\}$ coincides with the set of points where $A+B/B$ has rank $0$. Since the latter is locally free, that subset consists of certain connected components of $S$.
\end{rem}

\begin{lem}\label{lem:basistopology}
Let $(V,V^+)$ be as in the beginning of this section and let $T$ be an $S$-scheme.

Thus, the linear topology given by the submodules $\{A_T\vert A\sim V^+\}$ coincides with
the topology given by $\{{\mathcal A}\subset \hat V_T\vert {\mathcal A}\sim\hat V^+_T\}$.

In particular, if $W\to W'$ is a morphism of $\o_S$-modules continuous
w.r.t. the $V^+$-topology, then the induced morphism $\hat W_T\to \hat W'_T$ is also
continuous w.r.t. the $\hat V^+_T$-topology.
\end{lem}

\begin{pf}
Let $\mathcal A$ be a submodule of $\hat V_T$ such that ${\mathcal A}\sim\hat V^+_T$ and let $t\in T$ be a closed point. It suffices to show that there exists an open neighborhood $R$ of $t$ in $T$ and submodules $B,C\subset V$ commensurable with $V^+$ such that $B\vert_R\subseteq {\mathcal A}\vert_R\subseteq C\vert_R$. Indeed, since $V=\bigcup_{C\sim V^+}C$ it follows that
$V_{k(t)}=\bigcup_{C\sim V^+}C_{k(t)}$. Consider $C\sim V^+$ such
that ${\mathcal A}_{k(t)}\subset C_{k(t)}$. The previous proposition shows that this
inclusion also holds in a neighborhood $R$ of $t$ in $T$. The existence of
$B$ is proved similarly using the fact that $(0)=\bigcap_{B\sim
V^+}B$.
\end{pf}

\subsection{Grassmannian}

Let $T$ be an $S$-scheme, $L\subset \hat V_T$  a quasicoherent
$\o_T$-submodule, and $A\sim V^+$. We say that $(L,A)$ satisfies the
condition $(*)$ if
    $$(*)
    \hfill\hphantom{mmmmmmm}
    \hat V_T/L+\hat A_T\,=\,(0) \qquad \text{and}\qquad L\cap \hat A_T \text{ is l.f.f.t.}
    \hfill\hphantom{mmmmmmmm}
    $$
Observe that if $(L,A)$ satisfies $(*)$ ,then $L\cap A_T = L\cap \hat
A_T$ and $L=\hat L$.

\begin{lem}\label{lem:lema5partes}
Let $A,B\subseteq V$ be quasicoherent $\o_S$-submodules such that
$V/A, V/B$ are flat. Let $T$ be an $S$-scheme and $L\subseteq \hat
V_T$  a quasicoherent $\o_T$-submodule such that $\hat V_T/L$ is
flat. Assume that $(L,A)$ fulfills $(*)$. It then holds that
\begin{enumerate}
    \item $(L\cap \hat A_T)_{k(t)}^{\hat{\,}}=L_{k(t)}\cap \hat A_{k(t)}$ for all $t\in
    T$;
    \item if $A\subset B$, then $(L,B)$ satisfies $(*)$;
    \item if $\hat V_T/L+\hat B_T=0$, then $(L,B)$ satisfies $(*)$;
    \item $\hat V_T/L+\hat B_T$ is an $\o_T$-module of finite presentation;
    \item if $L\cap \hat B_T=0$, then $\hat V_T/L+\hat B_T$ es l.f.f.t.~.
\end{enumerate}
\end{lem}

\begin{pf}
Before proving the statements, let us note the following. Let
$C\subseteq D\subseteq V$ be arbitrary $\o_S$-submodules; then,
Snake's Lemma yields the following exact sequence
    \begin{equation}
    \label{eq:lemaser}
    \begin{aligned}
    0 \to L\cap\hat  C_T  & \to L\cap \hat D_T \to (D/C)_T \to \\
    & \to \hat V_T/L+\hat C_T \to \hat V_T/L+\hat D_T \to 0
    \end{aligned}
    \end{equation}
Let us proceed with the proofs.
\begin{enumerate}
    \item This follows from the flatness of $\hat V_T/L$ and the sequence
    $$
    0\to L\cap \hat A_T \to \hat A_T \to \hat V_T/L \to 0
    $$
    \item Sequence~(\ref{eq:lemaser}) implies that $\hat V_T/L+\hat B_T=(0)$
    and the exactness of
    $$
    0 \to L\cap \hat A_T   \to L\cap \hat B_T \to (B/A)_T \to  0
    $$
    Bearing in mind that $L\cap \hat A_T$ and $(B/A)_T$ are l.f.f.t.,  the claim follows.
    \item Point $(2)$ implies that  $(L,A+B)$ satisfies $(*)$. The sequence~(\ref{eq:lemaser}) for
    $B$ and $A+B$ reads
    $$
    0 \to L\cap \hat B_T   \to L\cap (A+B)_T^{\hat{\,}} \to (A+B/B)_T  \to 0
    $$
    Since the middle and the right hand side terms are l.f.f.t. one concludes.
    \item The sequence~(\ref{eq:lemaser}) for $B$ and $A+B$ gives
    $$
    L\cap (A+B)_T^{\hat{\,}} \to (A+B/B)_T  \to  \hat V_T/L+\hat B_T \to  0
    $$
    where the middle and the left hand side terms are l.f.f.t.~.
    \item The sequence~(\ref{eq:lemaser}) for $B$ and $A+B$ is
    $$
    0 \to L\cap (A+B)_T^{\hat{\,}}  \to (A+B/B)_T \to
    \hat V_T/L+\hat B_T \to 0
    $$
    where the middle and the left hand side terms are l.f.f.t.
    Further, the second part of this Lemma shows that $(L,A+B)$ satisfies $(*)$
    and, using the first part, one sees that
    {\footnotesize $$
    (L\cap (A+B)_T)_{k(t)}^{\hat{\,}} \,=\,  \hat L_{k(t)}\cap (A+B)_{k(t)}^{\hat{\,}}
    \,\hookrightarrow\,
    (A+B/B)_{k(t)} \,=\, (A_{k(t)}+B_{k(t)})/B_{k(t)}
    $$}
    or, what amounts to the same, $\hat V_T/L+\hat B_T$ is flat.
    These two facts allows one to prove that $\hat V_T/L+\hat B_T $ is
    l.f.f.t.~.
\end{enumerate}
\end{pf}

\begin{lem}\label{lem:sumadirecta}
Let $A\subseteq V$ be a quasicoherent $\o_S$-submodule such that
$V/A$ is flat and $A\sim V^+$. Let $T$ be an $S$-scheme and
$L\subseteq \hat V_T$  a quasicoherent $\o_T$-submodule such that
$\hat V_T/L$ is flat. Assume that $(L,A)$ satisfies $(*)$.

Thus, for every $t\in T$ there exists a neighborhood $R\subseteq T$
such that
    \begin{itemize}
    \item there exists $B\subset V$ such that $L_R\cap \hat B_R=0$ and $B\sim A$ is quasicoherent and $V/B$ flat;
    \item there exists $C\subset V$ such that $L_R\oplus \hat C_R=V_R$ and $C\sim A$ is quasicoherent and $V/C$ flat;
    \end{itemize}
\end{lem}

\begin{pf}
Let us see the first part. For any $B\subset A$ such that $A/B$ is
l.f.f.t. one has the exact sequence
    $$
    0 \to  L\cap \hat B_T \to  L\cap \hat A_T \to (A/B)_T
    $$
where the middle and right hand side terms  are l.f.f.t.~. Accordingly, let
$R$ be a neighborhood of $t$ such that $(L\cap \hat
A_T)_R^{\hat{\,}}=L_R\cap \hat A_R$ is a free $\o_R$-module of
finite rank. Comparing the above sequence for $B'\subset B$ with
$A/B', A/B$ l.f.f.t.  and recalling that $\bigcap_{B\sim V^+}
B=(0)$, the result follows.

Let us prove the second claim. Let $B$ and $R$ be as in the first
part. We therefore obtain the following exact sequence
    $$
    0 \to L_R\cap \hat A_R \to (A/B)_R \to
    \hat V_R/L_R+\hat B_R \to 0
    $$
where all terms are l.f.f.t. $\o_R$-modules. It follows that
$L_R\cap \hat A_R$ is an $R$-valued point of the Grassmannian of
$A/B$. Thus, shrinking $R$ if necessary, there exists a l.f.f.t.
$\o_S$-module, $\bar C\subseteq A/B$, such that $(A/B)/\bar C$ is
l.f.f.t. and
    \begin{equation}
    \label{eq:sumadirecta}
    (L_R\cap \hat A_R)\oplus \bar C_R\,=\, (A/B)_R
    \end{equation}

Let us check that $C:=\pi^{-1}(\bar C)\subseteq V$, where $\pi:A\to
A/B$, satisfies the requirements of the statement. Indeed,
equation~(\ref{eq:sumadirecta}) shows that the dashed arrow in the
commutative diagram
    {\footnotesize $$
    \xymatrix@C=14pt{
    0 \ar[r] & L_R\cap\hat C_R \ar[r] \ar[d] & L_R\cap \hat B_R =(0)\ar[r]
    \ar@{^(->}[d] & C_R /B_R \ar[r] \ar@{^(->}[d] \ar@{-->}[rd] & \hat V_R/L_R+\hat B_R \ar[r]
    \ar@{=}[d] & \hat V_R/L_R+\hat C_R \ar[r] \ar[d] & 0
    \\
    & 0 \ar[r] & L_R\cap \hat A_R \ar[r]
    & A_R /B_R \ar[r]  & \hat V_R/L_R+\hat B_R \ar[r]
    & 0
    }
    $$}
is an isomorphism. Hence, the kernel and cokernel of the dashed arrow vanish; i.e., $L_R\cap \hat C_R=(0)$ and $\hat
V_R/L_R+\hat C_R =(0)$.

Finally, since $A/C\simeq (A/B)/\bar C$ is l.f.f.t., we conclude by
Lemma~\ref{lem:lfftquotientimpliscomm}.
\end{pf}

\begin{lem}\label{lem:condicionequivGr}
Let $T$ be an $S$-scheme and $L\subseteq \hat V_T$ be a
quasicoherent $\o_T$-submodule such that $\hat V_T/L$ is flat.

If there exists $B\sim V^+$ with $V/B$ flat such that $L\cap \hat
B_T=(0)$ and $\hat V_T/L+\hat B_T$ is l.f.f.t., then for every point
$t$ there exists a neighborhood $R$ and a submodule $A\subseteq B$
such that $A\sim V^+$, $V/A$ is flat and $(L_R,\hat A_R)$ satisfies
$(*)$.
\end{lem}

\begin{pf}
Let us consider a submodule $A\subset V$ such that $A\sim V^+$, $V/A$ flat
and $B\subseteq A$. Let us consider the sequence
    {\small $$
    A/B \,\to \, V/B \,\to \, V/L+B
    $$}
Since $\cup_{A\sim V^+}=V$, there exists a neighborhood $R$ of $t$
and $A$ such that the above composition is surjective on $R$. We therefore have
    {\small $$
    0 \,\to \, L_R\cap \hat A_R \,\to \,
    (A/B)_R^{\hat{\,}} \,\to \,  \hat V_R/L_R+\hat B_R \,\to\, 0
    $$}
and, hence, $L_R\cap \hat A_R $ is l.f.f.t.~.
\end{pf}

\begin{thm}\label{thm:representGrV}
Let $V^+\subset V$ be sheaves of $\o_S$-modules as in the beginning
of this section. Accordingly, the functor from the category of $S$-schemes
to the category of sets defined by
    {\small $$
    \grv(T)\,:=\,
    \left\{\begin{gathered}
    \text{quasicoherent $\o_T$-submodules $L\subset \hat V_T$ s.t. $\hat V_T/L$ is flat and}
    \\
    \text{s.t. for each point of $T$ there exists a neighborhood $R$ and}
    \\
    \text{an $A\sim V^+$ with $\hat V_R/(L_R+\hat A_R)=0$ and $L_R\cap \hat A_R$ is l.f.f.t.}
    \end{gathered}
    \right\}
    $$}
is representable by an $S$-scheme that will be denoted by $\grv$ and will be
called the {\sl Grassmannian of $V$} (if $V^+$ is understood).
\end{thm}

\begin{pf}
From the very definition of the functor, one deduces that it is a
sheaf and that,  if can  therefore be assumed  that $S$ is an affine
scheme.  Let us consider the functors on the category of $S$-schemes
    {\small $$
    F_A(T)\,:=\, \left\{\begin{gathered}
    \text{quasicoherent $\o_T$-submodules}
    \\
    L\subset \hat V_T\text{ such that }L\oplus \hat A_T=\hat V_T
    \end{gathered}
    \right\}
    $$}
We shall prove that $F_A$ is representable for all $A\sim V^+$ and
that $\{F_A\}_{A\sim V^+}$ is a covering of $\grv$ by open
subfunctors.

Consider $B\sim V^+$ such that $B\subset A$. For each $C\sim V^+$
with $A\subset C$ the functor on the category of $S$-schemes defined
by
    {\small $$
    T\,\rightsquigarrow\,
    F_{C/A/B}(T):=\left\{\begin{gathered}
    \text{quasicoherent $\o_T$-modules  } L
    \\
    \text{ such that } L\oplus (A/B)_T \,=\, (C/B)_T
    \end{gathered}
    \right\}
    $$}
is representable by an $S$-scheme affine over $S$ (which is an open
subscheme of the Grassmannian of $C/B$). Observe that, as $C$
varies, these affine $S$-schemes form an inverse system whose
inverse limit represents the functor
    {\small $$
    T\,\rightsquigarrow\,
    F_{A/B}\,:=\,\left\{\begin{gathered}
    \text{quasicoherent $\o_T$-modules  } L
    \\
    \text{ such that } L\oplus (A/B)_T \,=\, (V/B)_T
    \end{gathered}
    \right\}\,=\,\limpl{C\supseteq A} F_{C/A/B}
    $$}
in the category of $S$-schemes. Let us denote this $S$-scheme by $F_{A/B}$
and let us observe that the map $F_{A/B}\to S$ is an affine morphism. Note
that as $B$ varies the schemes $\{F_{A/B}\}$ are an inverse system,
and since the maps to $S$ are affine, there exists its inverse
limit in the category of $S$-schemes. Finally, it is clear that
    \begin{equation}
    \label{eq:FAinverselimit}
    F_A\,=\, \limpl{B\subseteq A} F_{A/B}
    \end{equation}
For the sake of clarity, note that  there exists a correspondence
between $\hom(L,\hat A)$ and $F_A$ ($L$ being an $S$-valued point of
$F_A$), which assigns its graph to a map.

Secondly, Lemma~\ref{lem:sumadirecta} shows that the schemes
$\{F_A\}_{A\sim V^+}$ is a covering of $\grv$.

It remains to show that for a $S$-scheme $T$ and a functor
homomorphism $T\to \grv$ the map $F_A\times_{\grv} T \to T$ is
representable by an open subscheme of $T$ or, what amounts to the
same, that the set of points $t\in T$ s.t. $L_{k(t)}$ belongs to
$F_A(k(t))$ is open in $T$ (here $L$ is the module corresponding to
the map $T\to \grv$).

Let $t$ be such that  $L_{k(t)}\in F_A(k(t))$. Then, $\hat
V_{k(t)}/(L_{k(t)}+\hat A_{k(t)})=(0)$. From
Lemma~\ref{lem:lema5partes}, part 4, we know that $\hat V_T/L+\hat
A_T$ is of finite presentation and it therefore follows that $(\hat
V_T/L+\hat A_T)_t=(0)$. Let $R_0$ be the open set consisting of
those $t$ such that $(\hat V_T/L+\hat A_T)_t=(0)$. Thus, in $R_0$
the module $(L\cap  A)_{R_0}$ is l.f.f.t. and $L_t\cap\hat A_t=
(L\cap \hat A_T)_t$ for all $t\in R_0$. Therefore, the subset of those
points $t\in R_0$ such that $L_t\cap\hat A_t=(0)$ is the desired open
subscheme.
\end{pf}

\subsection{First Properties}

\begin{thm}
Let $T$ be an $S$-scheme. It then holds that $\gr(\hat
V_T)=\grv\times_S T$. In particular, for a closed point $s\in S$ one has that
$\gr(V)_{k(s)}=\gr(k(s)((t)))$.
\end{thm}

\begin{pf}
It is obvious that the canonical map $\grv\times_S T \hookrightarrow
\gr(\hat V_T)$ is an open immersion. It therefore suffices  to show that
$\gr(\hat V_{k(t)})=\grv\times_S \sp k(t)$ for a closed point $t\in
T$.

Let $L$ be a point in $\gr(\hat V_{k(t)})$ and ${\mathcal A}\subset
\hat V_{k(t)}$ with ${\mathcal A}\sim
\hat V_{k(t)}^+$ and $L\oplus {\mathcal A}= \hat V_{k(t)}$. Since
${\mathcal A}+\hat V_{k(t)}^+/ \hat V_{k(t)}^+$ is l.f.f.t. and
$\cup_{A\supseteq V^+} A=V$, there exists $A \sim V^+$ such that
${\mathcal A} \subseteq A_{k(t)}$. Now, the second item of
Lemma~\ref{lem:lema5partes} for ${\mathcal A}$ and $A_{k(t)}$
implies that $L$ is a $k(t)$-valued point of $\grv\times_S \sp
k(t)$.
\end{pf}

\begin{prop}\label{prop:GrVseparado}
The morphism $\grv\to S$ is separated.
\end{prop}

\begin{pf}
We follow \cite{EGAIHES}~I.\S5.5~. We must show that the diagonal
map $\grv\to \grv\times_S\grv$ is a closed immersion. This condition
is local, so it suffices to consider $A,B\sim V^+$ and an open subset $T\subset S$
such that $(A/A\cap B)_T, (B/A\cap B)_T$ are free of finite rank, and
to show that
    \begin{equation}\label{eq:separated}
    \o_{F_A}\otimes_{\o_T} \o_{F_B} \,\longrightarrow \, \o_{F_A\times_{\grv} F_B}
    \end{equation}
is surjective.

Let us describe $\o_{F_A\times_{\grv} F_B}$ in an explicit way. Let
$L\in\grv$ be the universal submodule. Thus, the canonical map
$L\oplus \hat A_{\grv}\to \hat V_{\grv}$ induces a map
    $$
    (A/A\cap B)_{\grv} \,\overset{\delta_{AB}}\longrightarrow\, \hat V_{\grv}/L+(A\cap B)_{\grv}^{\hat{\,}}
    $$
that, when restricted to $F_A$, is an isomorphism between l.f.f.t.
$\o_{F_A}$-modules of the same rank (since $A\sim B$ and part 5 of
Lemma~\ref{lem:lema5partes}).

Accordingly, on $F_A$ we consider the following composition of morphisms
    $$
    (B/A\cap B)_{F_A} \,\overset{\delta_{BA}}\longrightarrow\, \hat V_{F_A}/L_{F_A}+(A\cap B)_{F_A}^{\hat{\,}}
    \,\overset{\delta_{AB}^{-1}}\longrightarrow\, (A/A\cap B)_{F_A}
    $$
where all terms are l.f.f.t. $\o_{F_A}$-modules. Clearly, it holds
that the points of $F_A\times_{\grv} F_B$ are those points of $F_A$
where $ \delta_{AB}^{-1}\circ \delta_{BA}$ is an isomorphism. Hence,
$\o_{F_A\times_{\grv} F_B}$ is the localization of $\o_{F_A}$ by
$\det(\delta_{BA})$ or, what is tantamount, the localization of
$\o_{F_B}$ by $\det(\delta_{AB})$. And the map~(\ref{eq:separated})
is surjective.
\end{pf}

\begin{thm}\label{thm:imagendirectaOgrv}
Let $S$ be connected and let $\gr^0(V)$ denote a connected component
of $\grv$. Let $\pi\colon \gr^0(V)\to S$.

It therefore holds that
    \begin{enumerate}
    \item the canonical morphism $\o_S\to \pi_*\o_{\gr^0(V)}$ is an
isomorphism;
    \item $ H^0(\gr^0(V),\o_{\gr^0(V)})\,=\, H^0(S,\o_S)$;
    \item $ H^0(\gr^0(V),\o_{\gr^0(V)}^*)\,=\, H^0(S,\o_S^*)$.
    \end{enumerate}
\end{thm}

\begin{pf}
Recall that $\pi_*\o_{\gr^0(V)}$ is the sheaf of $\o_S$-modules associated to the presheaf
    {\small
    $$
    \, T\rightsquigarrow \,  H^0(\pi^{-1}(T),\o_{\gr^0(V)})\,=\, H^0(\gr^0(V_T),\o_{\gr^0(V_T)})
    $$}
where $T\subseteq S$ is an open subscheme. Therefore, if we show that
    $$
    H^0(\gr^0(V),\o_{\gr^0(V)})\,=\, H^0(S,\o_S)
    $$
holds for $S$ affine, then $(1)$ holds and therefore so does $(2)$.

Let us assume that $S$ is the spectrum of a local ring, $S=\sp\o_S$. Let $A\sim V^+$ be a
submodule such that $F_A\subset \gr^0(V)$. Let $L\in \gr^0(V)(S)$ be
a point of the Grassmannian such that $L\oplus \hat A=\hat V$; that
is, $L\in F_A(S)$. Let us fix elements $\omega\in\hom_{\o_S}(L,\o_S)$ and $a\in A$.

Now, the pair $(\omega,a)$ defines an affine line inside $F_A$ as
the image of
    $$
    \begin{aligned}
    j\colon {\mathbb A}^1_S & \hookrightarrow \,\hom_{\o_S}(L,A)\simeq F_A \\
    \lambda \, & \mapsto\quad  {\small (l\mapsto \lambda
    \omega(l)a )}
    \end{aligned}
    $$
Indeed, this procedure allows us to construct a projective line
inside the Grassmannian that, when intersected with $F_A$, coincides
with the previous construction. Indeed, bearing in mind that
$(L+<a>)/\ker\omega$ has rank $2$, we have
    $$
    \begin{aligned}
    \bar j\colon {\mathbb P}^1_S={\mathbb P}\big((L+<a>)/\ker\omega \big) & \hookrightarrow \gr^0(V) \\
    \bar L \qquad & \mapsto\,  p^{-1}(\bar L)
    \end{aligned}
    $$
where $p\colon (L+<a>)\to (L+<a>)/\ker\omega$. Note that $\ker\omega\in\grv$.

The restriction homomorphisms yield a commutative
diagram
    {\small $$
    \xymatrix{
    H^0(\gr^0(V),\o_{\gr^0(V)}) \ar[d]^{\bar j^*}\ar[r] &
    H^0(F_A,\o_{F_A})\simeq \o_S[\{x_i\}] \ar[d]^{ j^*}
    \\
    H^0({\mathbb P}^1_S,\o_{{\mathbb P}^1_S})= \o_S \ar[r] &
    H^0({\mathbb A}^1_S,\o_{{\mathbb A}^1_S})\simeq \o_S[y]
    }$$}
($\{x_i\}$ being a family of indeterminates).

Let a global section $s\in H^0(\gr^0(V),\o_{\gr^0(V)})$ be given. It
then holds that $(\bar j^*(s))\vert_{{\mathbb A}^1_S}=
j^*(s\vert_{F_A})$ and, in particular, $j^*(s\vert_{F_A})\in\o_S$;
i. e., it does not depend on $y$.

From expression~(\ref{eq:FAinverselimit}) and the very
definition of $j$, one sees that for each indeterminate $x_i$ of
the family $\{x_i\}$ there exists a pair $\omega,a$ such that
$j^*(x_i)$ depends on $y$. We conclude that $s\vert_{F_A}$ cannot
depend on $x_i$ for all $i$; that is, $s\vert_{F_A}\in\o_S$.

Since the above argument holds for all $A\sim V^+$ (with $F_A\subset
\gr^0(V)$) and since $\{F_A\}$ cover $\gr^0(V)$, it follows that
$s\in\o_S$. And the claim is proved.

For the third statement it suffices to show that
    $$
    H^0(\gr^0(V),\o_{\gr^0(V)}^*)\,=\, H^0(S,\o_S^*)
    $$
for any affine scheme $S$. This can be shown by applying similar arguments to those provided above to the diagram
    {\small $$
    \xymatrix{
    H^0(\gr^0(V),\o_{\gr^0(V)}^*) \ar[d]^{\bar j^*}\ar[r] &
    H^0(F_A,\o_{F_A}^*)\simeq (\o_S[\{x_i\}])^* \ar[d]^{ j^*}
    \\
    H^0({\mathbb P}^1_S,\o_{{\mathbb P}^1_S}^*)= \o_S^* \ar[r] &
    H^0({\mathbb A}^1_S,\o_{{\mathbb A}^1_S}^*)\simeq (\o_S[y])^*
    }$$}
\end{pf}

\section{The Determinant and Central Extensions}
\subsection{The Determinant}

\begin{thm}\label{thm:complexperfect}
Let $T$ be an $S$-scheme and let $L\in\grv(T)$ be a $T$-valued point
of $\grv$. Let $A\subseteq \hat V_T$ be a quasicoherent submodule
such that $\hat V_T/A$ is flat and $A\sim \hat V^+_T$.

Thus, the complex of $\o_T$-modules
    {\small $$
    \c^\bullet_A(L)\,\equiv\,\dots\to 0\to L\oplus A \overset{\delta_A}\longrightarrow \hat V_T\to 0
    \to \dots
    $$}
($\delta_A$ being the addition map) is a perfect complex (\cite{KM}).
\end{thm}

\begin{pf}
The claim is local in $T$ so it suffices to show that given $t$
there exists a neighborhood $R$ of $t$ in $T$ such that
$\c^\bullet_A(L)\vert_R= \c^\bullet_{\hat A_R}(L_R)$ is a perfect
complex. Given $t$, take $R$ such that $\hat A_R+\hat V^+_R/ \hat
V^+_R$ and $\hat A_R+\hat V^+_R/ \hat A_R$  are free of finite rank.
Observe that there is a quasiisomorphism of complexes (written
vertically)
    {\small $$
    \xymatrix{
    L_R\oplus \hat A_R \ar[r]^{p} \ar[d] & L_R \ar[d] \\
    \hat V_R \ar[r] & \hat V_R/\hat A_R
    }$$}
where $p$ denotes the projection.

Let us now, consider a quasicoherent submodule $B\subset V$ with $V/B$ flat
such that $B\sim V^+$ and $\hat A_R\subseteq \hat B_R$ (its existence
follows from Lemma~\ref{lem:basistopology}). Shrinking $R$ if necessary, we assume
that $\hat V_R/(L_R+\hat B_R)=(0)$ and that $L_R\cap \hat B_R$ is
l.f.f.t.~. Therefore, the exactness of the sequence~(\ref{eq:lemaser})
means that the morphism of complexes
    {\small $$
    \xymatrix{
    L_R\cap \hat B_R \ar[r] \ar[d] & L_R \ar[d] \\
    \hat B_R/\hat A_R  \ar[r] & \hat V_R/\hat A_R
    }$$}
is a quasiisomorphism. Since $\sim$ is an equivalence relation that
is compatible with pullbacks, we have that $\hat B_R\sim \hat
A_R$, and the result follows.
\end{pf}

Recalling the theory of determinant of complexes of \cite{KM}, we have the following

\begin{cor}\hfill{\,}
\begin{itemize}
    \item The determinant of the complex $\c^\bullet_A(L)$ changes base; that is, for a
    morphism of $S$-schemes $f\colon R\to T$ one has a canonical
    isomorphism
        {\small $$
        f^*\Det(\c^\bullet_A(L))\,\overset{\sim}\to\, \Det(\c^\bullet_{\hat A_R}(f^*L))
        $$}
    \item For $A,B\sim V^+_T$, there is a canonical isomorphism
        {\small $$
        \Det(\c^\bullet_A(L)) \otimes_{\o_T} \wedge (A/A\cap B)^*
        \,\overset{\sim}\to\,
        \Det(\c^\bullet_B(L)) \otimes_{\o_T} \wedge (B/A\cap B)^*
        $$}
        In particular, $\Det(\c^\bullet_A(L))$ defines a class on the relative Picard group
        that does not depend on the choice of $A$.
\end{itemize}
\end{cor}

From the second item of the previous corollary, we give the following

\begin{defn}
The determinant bundle on $\grv$ is the class of $\Det(\c^\bullet_A(\L))$ in the relative Picard group; that is, up to pullbacks of line bundles on $S$. A canonical representative of this class is $\Det(\c^\bullet_{V^+}(\L))$, which,  abusing  notation, will also be called the determinant bundle and will be denoted by $\Det_V$.
\end{defn}

\begin{rem}
Let $L$ be a $T$-valued point of $\grv$. The Euler-Poincar\'{e} characteristic of the complex $\c^\bullet_A(L)$
        {\small $$
        \begin{aligned}
        T\,& \longrightarrow\, \Z
        \\
        t \, &\longmapsto
        \chi(\c^\bullet_A(L)\otimes k(t))
        \end{aligned}
        $$}
which, for closed points is given by $\dim_{k(t)} L_{k(t)}\cap A_{k(t)}-\dim V_{k(t)}/(L_{k(t)}+A_{k(t)})$, is a locally constant function. Therefore, we have a map
        {\small $$
        \begin{aligned}
        \grv \,& \longrightarrow\, H^0(S,\Z)
        \\
        L \, &\longmapsto
        \chi(\c^\bullet_A(L))
        \end{aligned}
        $$}
and $H^0(S,\Z)$ can be used to label the connected components of $\grv$.
\end{rem}

\begin{defn}
Let $T$ be a $S$-scheme and $g:\grv_T\to \grv_T$ a homomorphism of
$T$-schemes. We say that $g$ preserves the determinant if there
exists a line bundle on $T$, which will be denoted by $\o_T(g)$, such that
    \begin{equation}\label{eq:isoGlenDet}
    g^*p_1^*\Det_V\,\simeq\,
    p_1^* \Det_V\otimes p_2^* \o_T(g)
    \end{equation}
on  $\grv_T:=\grv\times_S T$, where $\L$ is the universal submodule
and $p_i$ is the projection of $\grv\times_S T$ onto its $i$-th
factor.
\end{defn}

\begin{thm}\label{thm:elementpreserve}
Let $T$ be an $S$-scheme and let $g:\hat V_T\to\hat V_T$ be a
homomorphism of $\o_T$-modules such that
    \begin{itemize}
    \item $g$ is injective;
    \item $\hat V_T/g(\hat V^+_T)$ is flat and $\hat V^+_T \sim g(\hat V^+_T)$;
    \item $\hat V_T/g(\hat V_T)$ is flat and $\hat V_T \sim g(\hat V_T)$;
    \end{itemize}
Then, $g$ induces a  homomorphism of $T$-schemes
    $$
    g\colon  \grv_T\, \to\, \grv_T
    $$
which sends $L$ to $g(L)$ and preserves the determinant bundle.
\end{thm}

\begin{pf}
Let $L$ be a point of $\grv_T$. We must check that $g(L)$ also belongs to $\grv_T$. Observe that the left and right hand side terms of
the exact sequence
    {\small $$
    0\to \hat V_T/L \overset{g}\to \hat V_T/g(L)\to \hat V_T/g(\hat
    V_T)\to 0
    $$}
are quasicoherent flat $\o_T$-modules, the one in the middle also does.

On the other hand, since the condition $L\in\grv(T)$ is local in
$T$, we may assume that there exists $A\sim V^+$ such that $L\oplus
\hat A_T=\hat V_T$. Then, $g(L)\oplus g(\hat A_T)=g(\hat V_T)$ and,
therefore, $g(L)\cap g(\hat A_T)=(0)$ and
    {\small $$
    \hat V_T /g(L)+ g(\hat A_T)\,=\, \hat V_T / g(\hat V_T)
    $$}
is l.f.f.t. since $\hat V_T \sim g(\hat V_T)$. Recalling that
$g(\hat A_T)\sim g(\hat V^+_T)\sim \hat V^+_T$ and
Lemma~\ref{lem:condicionequivGr}, we conclude that $g(L)\in\grv_T$.

Let us now check that it preserves the determinant bundle. Let us now consider
 the following two exact sequences of complexes (written
vertically)
    {\small $$
    \xymatrix@C=13pt{
    0\ar[r] &
    \L_T\oplus \hat V^+_T \ar[d]\ar[r]^-g &
    g(\L_T)\oplus (\hat V^+_T + g(\hat V^+_T)) \ar[d] \ar[r] &
    (\hat V^+_T + g(\hat V^+_T))/g(\hat V^+_T) \ar[r] \ar[d] &
    0
    \\
    0 \ar[r] &
    \hat V_T \ar[r]^-g &
    \hat V_T \ar[r] &
    \hat V_T/ g(\hat V_T)\ar[r] &
    0}
    $$}
and
    {\small $$
    \xymatrix@C=13pt{
    0\ar[r] &
    g(\L_T)\oplus \hat V^+_T \ar[d]\ar[r] &
    g(\L_T)\oplus (\hat V^+_T + g(\hat V^+_T)) \ar[d] \ar[r] &
    (\hat V^+_T + g(\hat V^+_T))/\hat V^+_T \ar[r] \ar[d] &
    0
    \\
    0 \ar[r] &
    \hat V_T \ar[r] &
    \hat V_T \ar[r] &
    0 \ar[r] &
    0}
    $$}

Observe the following. The determinant of the complex on the
l.h.s. of the first sequence is the determinant bundle, $\Det_V$.
The complex on the r.h.s. of  the first sequence is perfect since
$\hat V^+_T \sim g(\hat V^+_T)$ and $\hat V_T \sim g(\hat V_T)$, and
its determinant is the pullback of a line bundle on $T$ by the projection
$\grv_T\to T$.

The complex in the middle of both sequences is perfect by
Theorem~\ref{thm:complexperfect}. The determinant of the complex on
the l.h.s. of the second sequence is $g^*\Det_V$, and the determinant
of the complex on the r.h.s. of the second sequence is the pullback
of a line bundle on $T$.

Now, the properties of the determinants give a canonical isomorphism
    {\footnotesize $$
    g^* p_1^*\Det_V\,\simeq\, p_1^*\Det_V\otimes
    \wedge \big(\hat V^+_T + g(\hat V^+_T)/\hat V^+_T\big)^*\otimes
    \wedge \big(\hat V^+_T + g(\hat V^+_T)/g(\hat V^+_T)\big)\otimes
    \wedge \big(\hat V_T /g(\hat V_T)\big)
    $$}
(where $\wedge$ denotes the exterior algebra of highest degree) and
the result follows.
\end{pf}

\begin{rem}\label{rem:conditonsThmOpen}
Let $T$ be an $S$-scheme and let $g:\hat V_T\to\hat V_T$ be an injective
homomorphism of $\o_T$-modules. Therefore, the set of points
where the second and third conditions are satisfied is an open
subscheme.
\end{rem}

\subsection{The Extension}\label{subsect:extension}

\begin{defn}
Let $\grv\to S$ be the Grassmannian of $V$. For each $S$-scheme $T$, we define the sets
        {\small $$
        {\mathcal Q}_S(T)
        \,:=\,
        \left\{ \begin{gathered}
        \text{$\o_T$-module monomorphisms $\hat V_T\to\hat V_T$}
        \\
        \text{such that the conditions of theorem~\ref{thm:elementpreserve} hold}
        \end{gathered}\right\}
        $$}
and
        {\small $$
        \widetilde{\mathcal Q}_S(T)
        \,:=\,
        \left\{
        \begin{gathered}\text{pairs $(\tilde g, g)$ where }g\in {\mathcal Q}(T)
        \\
        \text{ and }\tilde g\in\aut_{\o_T} \o_T(g)
        \end{gathered}
        \right\}
        $$}
where
    {\small
    $$\o_T(g):= \wedge \big(\hat V^+_T + g(\hat V^+_T)/\hat V^+_T\big)^*\otimes
    \wedge \big(\hat V^+_T + g(\hat V^+_T)/g(\hat V^+_T)\big)\otimes
    \wedge \big(\hat V_T /g(\hat V_T)\big)\, .$$}
\end{defn}

\begin{thm}
With the above notation, it holds that
    \begin{enumerate}
        \item ${\mathcal Q}_S(T)$ is a monoid.
        \item The assignment $T\rightsquigarrow {\mathcal Q}_S(T)$ is a contravariant functor
        from the category of $S$-schemes to the category
        of monoids,
        \item The map $T\rightsquigarrow {\mathcal Q}_S(T)$, where $T$ is an open subset of $S$,
        defines a sheaf on $S$, which will be denoted by ${\mathcal Q}_S$.
    \end{enumerate}
The same properties hold for $\widetilde{\mathcal Q}_S$.
\end{thm}

\begin{pf}
$(1)$ Observe that the identity $\hat V_T\to \hat V_T$ belongs to ${\mathcal Q}_S(T)$ and that
if $g,h\in {\mathcal Q}_S(T)$ satisfy the conditions of the Theorem, then
so does $g\circ h$. The only non-trivial part is that $\hat
V_T/(g\circ h)(\hat V_T)$ is flat. Note that the left and right hand
side terms of the exact sequence
    $$
    0\to \hat V_T/h(\hat V_T)\,\overset{g}\to\,
    \hat V_T/(g\circ h)(\hat V_T)\,\to\,
    \hat V_T/g(\hat V_T)\,\to\, 0
    $$
are flat. So is the one in the middle.

Let us detail explicitly the composition law of $\widetilde{\mathcal Q}_S$. Let $(f,\tilde f)$ and $(g,\tilde g)$ be two elements of $\widetilde{\mathcal Q}_S(T)$. Since $\tilde g\in\aut_{\o_S}  \o_T(g) $, it follows that $f^*(\tilde g)\in \aut_{\o_T}  f^*\o_T(g)$. Note that
    {\small
    $$
    \begin{aligned}
    g^*\o_T(f) \otimes \o_T(g) \, & \simeq\,
    g^* \big( f^* p_1^*\Det_V\otimes (p_1^*\Det_V)^{-1}\big) \otimes g^* p_1^*\Det_V\otimes (p_1^*\Det_V)^{-1}
    \,\simeq \\
    & \simeq \,
    (gf)^* p_1^*\Det_V \otimes (p_1^*\Det_V)^{-1}
    \end{aligned}
    $$}
and  therefore $g^*(\tilde f)\otimes \tilde g\in \aut_{\o_T}(\o_T(gf))$. Summing up, the composition law is
    $$
    (g,\tilde g) (f,\tilde f)
    \,:=\,
    (gf, g^*(\tilde f)\otimes \tilde g)
    $$

$(2)$  Let $T$ be an $S$-scheme and
let $g$ be an element in ${\mathcal Q}_S(T)$. If $R\to T$ is a morphism
of $S$-schemes, then the morphism $g\otimes 1:\hat
V_T\hat\otimes_{\o_T}\o_R\to\hat V_T\hat\otimes_{\o_T}\o_R$
satisfies the conditions of Theorem~\ref{thm:elementpreserve}. And the homomorphism of $R$-schemes
induced by $g\otimes 1:\grv_R\to \grv_R$ coincides with the one induced
by $g:\grv_T\to \grv_T$ by $R\to T$. This means that there is a
canonical map ${\mathcal Q}_S(T)\to {\mathcal Q}_S(R)$, which
is an homomorphism of monoids.

$(3)$ Straightforward.
\end{pf}

We shall now study the map $\widetilde{\mathcal Q}_S\to {\mathcal Q}_S$ from the point of view of Leech's theory of {\sl extensions of groups by moniods} (see~\cite{Leech}). We shall follow his notations and results. Recall that a {\sl normal extension} of the group $K$ by the monoid $Q$ is a sequence $K\overset{i}\to\tilde Q\overset{\pi}\to Q$, where $\tilde Q$ is a monoid, $K$ is a subgroup of the group of units of $\tilde Q$ such that $qK=Kq$ for all $q\in\tilde Q$, $i$ is an injective morphism of monoids, and $p$ is a surjective morphism of monoids such that it induces an isomorphism of the quotient monoid of left cosets $\tilde Q/K$ with $Q$.

\begin{thm}
For an $S$-scheme $T$, there is a normal extension of monoids
    \begin{equation}\label{eq:normalextension}
    H^0(T,\o_T^*) \,\longrightarrow \, \widetilde{\mathcal Q}_S(T) \,\overset{\pi_T}\longrightarrow\, {\mathcal Q}_S(T)
    \end{equation}
\end{thm}

\begin{pf}
By the very definition the map $\widetilde{\mathcal Q}_T \to {\mathcal Q}_T$ is surjective morphism of monoids.
The fiber of any element of ${\mathcal Q}_T $ is isomorphic to $H^0(T,\o_T^*)$ by Theorems~\ref{thm:elementpreserve} and~\ref{thm:imagendirectaOgrv}.
\end{pf}

In particular, one could understand the above normal extension as a normal extension of the sheaf on monoids ${\mathcal Q}_S$ by the sheaf on groups ${\mathbb G}_{m,S}$ (the relative multiplicative group).

\begin{rem}
The above  definition of $\widetilde{\mathcal Q}_S$ is closer to that employed in \cite{AP} than to that of \cite{MPa}. However, let us see that  essentially both coincide.  Let us denote by $\overline{\mathcal Q}_S$ the extension of the latter reference, which consists of pairs $(g,\bar g)$ where $\bar g$ is an isomorphism $g^* p_1^*\Det_V\simeq (p_1^*\Det_V)\otimes \o_T(g)$. Let us fix any set-theoretic section of $\overline{\mathcal Q}_S \to {\mathcal Q}_S$, say $\sigma$. Accordingly,  the map
    {\small $$
    \begin{aligned}
    \widetilde{\mathcal Q}_S \, & \longrightarrow \,\overline{\mathcal Q}_S
    \\
    (g,\tilde g) \, & \longmapsto \, (g, (1\otimes \tilde g)\circ \sigma(g))
    \end{aligned}
    $$}
can be used to identify both central extensions as follows
    {\small $$
    \xymatrix@R=5pt{
    & & \widetilde{\mathcal Q}_S \ar[dd]^{\wr} \ar[dr]
    \\
    0 \ar[r] & {\mathbb G}_{m,S} \ar[ru] \ar[dr] & & {\mathcal Q}_S \ar[r] & 0
    \\
    & & \overline{\mathcal Q}_S \ar[ur]
    }
    $$}
\end{rem}

Following Leech's paper, we know that the extension
    $$
    H^0(T,\o_T^*) \,\longrightarrow \, \widetilde{\mathcal Q}_T \,\overset{\pi_T}\longrightarrow\, {\mathcal Q}_T
    $$
arises from two functors $F_T$ (from the category given by the ${\mathcal L}$-quasiorder on $Q$ to the category of groups) and $G_T$ (from the category given by the ${\mathcal R}$-quasiorder on $Q$ to the category of groups) and a factor system
    $$
    \alpha_T\,\colon\, {\mathcal Q}_T \times {\mathcal Q}_T \,\longrightarrow\,
    \bigcup_{x\in {\mathcal Q}_T} F_T(x)
    $$
Therefore, the extension $\widetilde{\mathcal Q}_T$ is the set
$\bigcup_{x\in {\mathcal Q}_T} \{x\}\times F_T(x)$ endowed with the following composition law
    $$
    (x,a)\cdot (y,b)\,:=\,
    (xy, G_x^{xy}(a) \alpha_T(x,y) F_y^{xy}(b))
    $$
where $G_x^{xy}$ denotes the image of the morphism $x \geq_{\mathcal R} xy$ by the functor $G_T$; i.e. $G_x^{xy}: G_T(x)\to G_T(xy)$, and analogously for $F_y^{xy}$ since $y \geq_{\mathcal L} xy$.

From the very construction of $\widetilde{\mathcal Q}_T$ we deduce the following properties
    \begin{itemize}
    \item since $\widetilde{\mathcal Q}_T$ is associative, it holds that $\alpha_T(1,x)=\alpha_T(x,1)=1_{F_T(1)}$.
    \item $F_T(x)=\pi^{-1}(x)= H^0(T,\o_T^*)$ (Theorems~\ref{thm:elementpreserve} and~\ref{thm:imagendirectaOgrv}); that is $F_T$ is the constant functor  $x\rightsquigarrow H^0(T,\o_T^*)$, which is a commutative group (same for $G_T$).
    \end{itemize}
In light of this latter observation, and by Theorem 2.13 of~\cite{Leech}, we know that our extension {\sl splits} (i.e. there exists a section of $\pi_T$ that is a morphism of monoids) if and only if there exists $\phi\in \prod_{x\in{\mathcal Q}_T} F_T(x)$ with $\phi(1)=1$, such that
    $$
    \alpha_T(x,y)\,=\,
    (\delta\phi)(x,y)\,:=\, G_x^{xy}(\phi(x))\phi(xy)^{-1} F_y^{xy}(\phi(y))
    $$

Let $f:T\to S$ be given. We can therefore consider the transformation of the extension $H^0(S,\o_S^*) \to \widetilde{\mathcal Q}_S \to {\mathcal Q}_S$ by $H^0(S,\o_S^*)\to H^0(T,\o_T^*)$. This extension, which we denote by $f^* \widetilde{\mathcal Q}_S$, is that associated with the constant functors $F=G=H^0(T,\o_T^*)$ and the factor system
    $$
    f^*\alpha_S: {\mathcal Q}_S \times {\mathcal Q}_S \,\overset{\alpha_S}\longrightarrow\,
    \bigcup_{x\in {\mathcal Q}_S} H^0(S,\o_S^*)\,\longrightarrow\, \bigcup_{x\in {\mathcal Q}_S} H^0(T,\o_T^*)
    $$
and is therefore of the following type
    $$
    H^0(T,\o_T^*) \,\to \, f^* \widetilde{\mathcal Q}_S \,\to\, {\mathcal Q}_S
    $$

On the other hand, we can transform $H^0(T,\o_T^*) \to \widetilde{\mathcal Q}_T \to {\mathcal Q}_T$ by $f^*:{\mathcal Q}_S \to {\mathcal Q}_T $. This extension, which we denote by $f_* \widetilde{\mathcal Q}_T $, is that associated with the constant functors $F=G=H^0(T,\o_T^*)$ and the factor system
    $$
    \begin{aligned}
    f_*\alpha_T: {\mathcal Q}_S  \times {\mathcal Q}_S \,&  \longrightarrow\,
    \bigcup_{x\in {\mathcal Q}_S} H^0(T,\o_T^*)
    \\
    (x,y) \, &\longmapsto \, \alpha_T(f^*(x),f^*(y))
    \end{aligned}$$
and it has the form
    $$
    H^0(T,\o_T^*) \,\to \, f_* \widetilde{\mathcal Q}_T \,\to\, {\mathcal Q}_S
    $$

The explicit expressions show that both extensions coincide; that is, $f^* \widetilde{\mathcal Q}_S =
f_* \widetilde{\mathcal Q}_T$, or, what is tantamount, the following diagram is commutative
    $$
    \xymatrix{
    {\mathcal Q}_S \times {\mathcal Q}_S \ar[d]_{f^*\times f^*} \ar[r]^-{\alpha_S}
    &
    \bigcup_{x\in {\mathcal Q}_S} H^0(S,\o_S^*) \ar[d]
    \\
    {\mathcal Q}_T \times {\mathcal Q}_T  \ar[r]^-{\alpha_T}
    &
    \bigcup_{x\in {\mathcal Q}_T} H^0(T,\o_T^*)
    }$$

Now, if we consider $T$ varying in the open subschemes of $S$, we have proved the following

\begin{thm}\label{thm:factorsystemQsheaf}
The factor system defines a sheaf homomorphism
    $$
    {\mathcal Q}_S \times {\mathcal Q}_S \,\overset{\alpha}\longrightarrow \,{\mathbb G}_{m,S}
    $$
\end{thm}


Inspired by \cite{MPa} we give the following criterion for splitting.

\begin{thm}
Let $T$ be an $S$-scheme and ${\mathcal A}\in \grv(T)$ be a point. Let ${\mathcal Q}_T^{\mathcal A}$ be the submonoid of ${\mathcal Q}_T$ consisting of those $g$ such that $g({\mathcal A})= {\mathcal A}$.

Thus, the normal extension of ${\mathcal Q}_T^{\mathcal A}$ induced by~(\ref{eq:normalextension})
    $$
        H^0(T,\o_T^*) \,\overset{\iota}\to \, \widetilde{\mathcal Q}_T^{\mathcal A} \,\to\, {\mathcal Q}_T^{\mathcal A}
    $$
splits.
\end{thm}

\begin{pf}
First, let us note that the inclusion ${\mathcal Q}_T^{\mathcal A} \hookrightarrow {\mathcal Q}_T$ gives rise to the above normal extension.

Since $\widetilde{\mathcal Q}_T^{\mathcal A}$ consists of pairs $(g,\tilde g)$ such that $g({\mathcal A})={\mathcal A}$, it follows that it acts on $\Det_V\times_{\grv_T}\{{\mathcal A}\}$, which is a locally free $\o_T$-module of rank 1. Thus, we obtain a homomorphism of monoids
    $$
    \tilde\mu\colon \widetilde{\mathcal Q}_T^{\mathcal A} \,\longrightarrow\, H^0(T,\o_T^*)
    $$

The very construction shows that the composition
    $$
    H^0(T,\o_T^*) \,\overset{\iota}\to \, \widetilde{\mathcal Q}_T^{\mathcal A} \,\overset{\tilde\mu}\to\, H^0(T,\o_T^*)
    $$
is the identity.

Let $\sigma$ be a section of $\pi$ as sets and let $\phi$ be defined by the relation $\sigma(x)=(x,\phi(x))\in \{x\}\times F_T(x)$ for $x\in {\mathcal Q}_T$. We normalize $\sigma$ by the condition $\phi(1)=1$. Let $c$ denote the cocycle associated with $\sigma$ which is the only element $c(x,y)\in H^0(T,\o_T^*)$ satisfying the following identity
    $$
    \sigma(x)\sigma(y)\, =\,
    c(x,y)\sigma(x y)
    $$
We claim that $G_1^{xy}(c(x,y))= (\alpha+\delta\phi)(x,y)$, where the symbol $+$ on the r.h.s. is  the addition of factor systems, which is defined pointwise. Indeed, this  follows from the equality of the following two identities in $\tilde{\mathcal Q}_T$
    {\small $$
    \begin{gathered}
    \sigma(x)\sigma(y)\, =\,
    (x,\phi(x))\cdot (y,\phi(y)) \,=\,
    \big(xy, G_x^{xy}(\phi(x))\alpha_T(x,y)F_y^{xy}(\phi(y))\big)
    \\
    c(x,y)\sigma(x y)\,=\,
    \big(xy,G_1^{xy}(c(x,y))\phi(xy)\big)
    \end{gathered}
    $$}

Consider the following section of $\pi$
    $$
    \sigma'(x)\,:=\, \tilde\mu(\sigma(x))^{-1} \cdot\sigma(x)
    $$
Let us compute the cocycle, $c'$, associated with $\sigma'$. Recall that $c'$ satisfies
    $$
    \sigma'(x)\sigma'(y)\, =\,
    c'(x,y)\sigma'(x y)
    $$
and applying $\tilde \mu$ we obtain
    $$
    \tilde\mu(\sigma(x))^{-1} \sigma(x)
    \tilde\mu(\sigma(y))^{-1} \sigma(y)
    \,=\,
    c'(x,y)
    \tilde\mu(\sigma(x y))^{-1} \sigma(x y)
    $$
and therefore
    $$
    c(x,y)  \sigma(x y)
    \tilde\mu(c(x,y))^{-1}\tilde\mu(\sigma(x y))^{-1}
    \,=\,
    c'(x,y)
    \tilde\mu(\sigma(x y))^{-1} \sigma(x y)
    $$
Simplifying, we have
    $$
    c(x,y)
    \tilde\mu(c(x,y))^{-1}
    \,=\,
    c'(x,y)
    $$
Recalling that $\tilde\mu$ is the identity on the elements of $H^0(T,\o_T^*)$, it follows that $c'\equiv 1$; or, what is tantamount, that the section $\sigma'$ is a morphism of monoids. The statement is proved.
\end{pf}

\subsection{Enlargement}

Note that the conditions for a monomorphism to belong to ${\mathcal Q}_S$ are rather strong and that the origin of these conditions is Theorem~\ref{thm:elementpreserve}. Motivated by Remark~\ref{rem:conditonsThmOpen}, we wonder whether they could be weakened; that is, if the monoid ${\mathcal Q}_S$ could be enlarged. Indeed, we shall show that there are cases in which this can be done.

Henceforth, $S$ is a non-singular curve over an algebraically closed field (not necessarily complete). Then, Ex 6.11, Chp. II of \cite{Hart} shows that  there is an isomorphism
    $$
    (\Det,\rk)\colon K(S)\,
    \overset{\sim}\longrightarrow\,
    \operatorname{Pic}(S)\times {\mathbb Z}
    $$
where $K(S)$ denotes the Grothendieck group of coherent sheaves of $\o_S$-modules. Furthermore, this is also valid for any non-empty open subset of $S$.

\begin{rem}
Another case in which this approach holds (because of the structure of its Grothendieck group) is the case of Dedekind domains. In a future work we plan to study this case and its applications to arithmetic (e.g. towards a reciprocity law for families of curves generalizing the ideas of \cite{AP,MPa})
\end{rem}

Let $[M]\in K(S)$ denote the class of a coherent sheaf $M$ in $K(S)$. Let us consider the category whose objects are complexes
    $$
    M_{\bullet}\equiv \ldots\to 0 \to M_0 \to M_1 \to 0 \to\ldots
    $$
where $M_i$ is $i$-th degree part of $M_{\bullet}$; $M_i$ is a quasicoherent $\o_S$-submodule of $V$; $M_0\subseteq M_1$ and the map $M_0\to M_1$ is the inclusion; and the cokernel, $M_1/M_0$, is coherent. The morphisms are the standard morphisms of complexes of $\o_S$-modules.

To each object $M_{\bullet}$ on this category we may attach an element of the Grothendieck group, namely, the class $[M_1/M_0]\in K(S)$. Note that an exact sequence $0\to M'_{\bullet}\to M_{\bullet}\to M''_{\bullet}\to 0$ induces an exact sequence between the cokernels and, hence, yields the identity  $[M_1/M_0]=[M'_1/M'_0]+[M''_1/M''_0]$ in $K(S)$. Similarly, if $M_{\bullet}\to N_{\bullet}$ is a quasi-isomorphism, then $M_{\bullet}$ and $N_{\bullet}$ give rise to the same class $[M_1/M_0]= [N_1/N_0]$.

Let us introduce the sheaf on $S$
        {\small $$
        \bar{\mathcal Q}_S(T)
        \,:=\,
        \left\{ \begin{gathered}
        \text{$\o_T$-module monomorphisms $g:\hat V_T\to\hat V_T$ such that }
        \\
        \hat V_T/ g \hat V_T^+ \text{ is quasicoherent and } (\hat V_T^++g \hat V_T^+)/\hat V_T^+ ,
        \\
        (\hat V_T^++g \hat V_T^+)/g \hat V_T^+ ,
        \hat V_T/ g \hat V_T
        \text{ are coherent}
        \end{gathered}\right\}
        $$}
for $T\subseteq S$, an open subscheme.

Observe that this monoid  contains the {\sl restricted general linear group} of \cite{AP} and that it could therefore be used to reformulate their results in our setup.

\begin{thm}
$\bar{\mathcal Q}_S(T)$ is a sheaf on monoids.
\end{thm}

\begin{pf}
We must show that given $f,g\in \bar{\mathcal Q}_S(T)$, then $g
f\in \bar{\mathcal Q}_S(T)$. First, note  that all sheaves
involved are quasicoherent since $V$ and $V^+$ are too and that it
therefore suffices to show that certain (quasicoherent)
sheaves are coherent.

Note that the hypotheses imply that the r.h.s. and l.h.s. terms of
the exact sequence
    {\small
    $$\xymatrix@=15pt{
    0 \ar[r] & \hat V_T / f\hat V_T \ar[r]^g &  \hat V_T / gf\hat V_T \ar[r] &
    \hat V_T / g\hat V_T \ar[r] & 0}
    $$}
are quasicoherent and, hence, we have that $\hat V_T / gf\hat
V_T$ is also coherent.

Let us show that $(\hat V_T^++g f \hat V_T^+)/\hat V_T^+ $  is
coherent. The module $(\hat V^+ + g\hat V^+ + gf \hat V^+)/(\hat
V^++ g \hat V^+)$ is coherent since there is a surjection
    {\small
    $$\xymatrix@=18pt{
    (\hat V^+ + f \hat V^+)/ \hat V^+   \ar@{->>}[r]^-{g} &
    (\hat V^+ + g\hat V^+ + gf \hat V^+)/(\hat V^++ g \hat V^+)}
    $$}
Thus, from the exact sequence
    {\scriptsize
    $$\xymatrix@=10pt{
    0 \ar[r] & (\hat V^+ + g \hat V^+)/ \hat V^+  \ar[r] &   (\hat V^+ + g\hat V^+ + gf \hat V^+)/\hat V^+ \ar[r] &
    (\hat V^+ + g\hat V^+ + gf \hat V^+)/(\hat V^++ g \hat V^+) \ar[r] & 0}
    $$}
we have that the middle term is also coherent since both sides are too.
Note  that $(\hat V^+ + f\hat V^+ + gf \hat V^+)/(\hat V^++ g f
\hat V^+)$ is coherent since
    {\scriptsize
    $$\xymatrix@=18pt{
    (\hat V^+ + f \hat V^+)/ f \hat V^+   \ar@{->>}[r]^-{g} &
    (\hat V^+ + g\hat V^+ + gf \hat V^+)/(\hat V^++ g f  \hat V^+)}
    $$}
The latter two facts says that the middle and r.h.s. terms of the sequence
    {\scriptsize
    $$\xymatrix@=10pt{
    0 \ar[r] & (\hat V^+ + g f \hat V^+)/ \hat V^+  \ar[r] &   (\hat V^+ + g\hat V^+ + gf \hat V^+)/\hat V^+ \ar[r] &
    (\hat V^+ + g\hat V^+ + gf \hat V^+)/(\hat V^++ g f \hat V^+) \ar[r] & 0}
    $$}
are coherent, and hence so is the l.h.s. term; that is, $(\hat V^+ + g f \hat V^+)/ \hat V^+$ is coherent.

It remains to show that $(\hat V_T^++g f \hat V_T^+)/g f \hat
V_T^+ $  is coherent. The arguments are similar to those given
above and use the following facts. The module $(\hat V^+ + g \hat
V^+ + gf \hat V^+)/(g \hat V^++ g f \hat V^+)$ is coherent since
it is a quotient of $(\hat V^+ +  g \hat V^+)/ g \hat V^+$. The
module $(\hat V^+ + g\hat V^+ + gf \hat V^+)/( \hat V^++ g f \hat
V^+)$ is coherent since it is a quotient of $(\hat V^+ +  g \hat
V^+)/ \hat V^+$. There are exact sequences
    {\scriptsize
    $$\xymatrix@=13pt{
    0 \ar[r] & (\hat V^+ + f \hat V^+)/ f\hat V^+ \ar[r]^-g &
    (\hat V^+ + g \hat V^+ + g f \hat V^+)/g f \hat V^+ \ar[r] &
    (\hat V^+ + g \hat V^+ + g f \hat V^+)/( g \hat V^++ g f \hat V^+) \ar[r] &
    0}
    $$}
and
    {\scriptsize
    $$\xymatrix@=13pt{
    0 \ar[r] & (\hat V^+ + g f \hat V^+)/g f\hat V^+  \ar[r] &
    (\hat V^+ + g \hat V^+ + g f \hat V^+)/g f \hat V^+ \ar[r] &
    (\hat V^+ + g \hat V^+ + g f \hat V^+)/(\hat V^++ g f \hat V^+) \ar[r] & 0}
    $$}
\end{pf}

\begin{defn}
We associate the element of $K(T)$ defined by
    $$
    [f]\,:=\,
    -\left[ (\hat V_T^++f \hat V_T^+)/\hat V_T^+ \right] +
    \left[ (\hat V_T^++f \hat V_T^+)/f \hat V_T^+ \right] -
    \left[ \hat V_T/ f \hat V_T\right]
    \,\in\, K(T)
    $$
to an element $f\in \bar{\mathcal Q}_S(T)$.
\end{defn}

Given an element  $g\in \bar{\mathcal Q}_S$ and an object
$M_{\bullet}$, we may consider the object $g_* M_{\bullet} :=
gM_0\to g M_1$. Therefore, the map $M_{\bullet} \to g_* M_{\bullet}$
gives rise, by linearity, to a transformation defined on the elements of the
above type
    $$
    [f]\,  \rightsquigarrow \,g_* [f]
    $$

\begin{thm}
Let $g, f\in \bar{\mathcal Q}_S(T)$. It holds that
    $$
    g_*[f] \,=\, [gf] - [g]
    $$
\end{thm}

\begin{pf}
Let $A\subseteq \hat V_T$ be a submodule such that $\hat V_T^+,
f\hat V_T^+, g\hat V_T^+, gf\hat V_T^+ \subseteq A$ and such that
the middle terms of sequences $(2)$ to $(7)$ below are coherent.
Note that the previous result shows that $A=\hat V_T^+ + f\hat
V_T^+ + g\hat V_T^+ + gf\hat V_T^+$ does the job.

Let us consider the following set of exact sequences
    {\small
    $$\xymatrix@=13pt{
    0 \ar[r] & \hat V_T / f\hat V_T \ar[r]^-g &  \hat V_T / gf\hat V_T \ar[r] &
    \hat V_T / g\hat V_T \ar[r] & 0 & (1)
    \\
    0 \ar[r] & \hat V_T^++f\hat V_T^+ / f\hat V_T^+ \ar[r]^-g &  A / gf\hat V_T^+ \ar[r] &
    A/ (g\hat V_T^+ + gf\hat V_T^+ ) \ar[r] & 0  & (2)
    \\
    0 \ar[r] & \hat V_T^++gf\hat V_T^+ / gf\hat V_T^+ \ar[r] &  A / gf\hat V_T^+ \ar[r] &
    A/ (\hat V_T^+ + gf\hat V_T^+) \ar[r] & 0 & (3)
    \\
    0 \ar[r] & \hat V_T^++g\hat V_T^+ / g\hat V_T^+ \ar[r] &  A / g\hat V_T^+ \ar[r] &
    A/ (\hat V_T^+ + g\hat V_T^+ ) \ar[r] & 0 & (4)
    \\
    0 \ar[r] & \hat V_T^++f\hat V_T^+ / \hat V_T^+ \ar[r]^-g &  A / g\hat V_T^+ \ar[r] &
    A/ (g\hat V_T^+ + gf\hat V_T^+ ) \ar[r] & 0 &(5)
    \\
    0 \ar[r] & \hat V_T^++gf\hat V_T^+ / \hat V_T^+ \ar[r] &  A / \hat V_T^+ \ar[r] &
    A/ (\hat V_T^+ + gf\hat V_T^+ ) \ar[r] & 0 &(6)
    \\
    0 \ar[r] & \hat V_T^++g\hat V_T^+ / \hat V_T^+ \ar[r] &  A / \hat V_T^+ \ar[r] &
    A/ (\hat V_T^+ + g\hat V_T^+ ) \ar[r] & 0 &(7)
    }$$
    }
Observe, furthermore, that all sheaves are quasicoherent and that
by hypothesis and by the previous results all terms on left hand
sides are coherent and all middle terms are coherent sheaves by the
choice of $A$. Thus, all terms on the right hand sides are coherent
too.

The class $g_*[f]$ is, by definition of $g_*$, equal to
    {\small $$
    g_*[f] \,=\,
    -\left[ (g\hat V_T^++gf \hat V_T^+)/g\hat V_T^+ \right] +
    \left[ (g\hat V_T^++gf \hat V_T^+)/gf \hat V_T^+ \right] -
    \left[ g\hat V_T/ gf \hat V_T\right]
    $$}
Accordingly, computing the first term with the help of the
sequence $(5)$, the second one with the sequence $(2)$ and the
last one with the sequence $(1)$, one has that
    {\small
    $$
    \begin{aligned}
    g_*[f] \,& =\,
    -[A / g\hat V_T^+]+[A/ (g\hat V_T^+ + gf\hat V_T^+ )] +
    [A / gf\hat V_T^+] - [A/ (g\hat V_T^+ + gf\hat V_T^+ )] -
    \\
    &\qquad\qquad - [ \hat V_T / gf\hat V_T] + [\hat V_T / g\hat V_T]
    \,=
    \\
    & = \,
    -[A / g\hat V_T^+]+
    [A / gf\hat V_T^+] -
    [ \hat V_T / gf\hat V_T] + [\hat V_T / g\hat V_T]
    \end{aligned}
    $$}
Now, plug in the relations given by sequences $(4)$ and $(3)$ into
the last expression
    {\small
    $$
    \begin{aligned}
    g_*[f] \,& =\,
    - [\hat V_T^++g\hat V_T^+ / g\hat V_T^+] - [A/ (\hat V_T^+ + g\hat V_T^+ )]
    + [\hat V_T^++gf\hat V_T^+ / gf\hat V_T^+] +
     \\
    &\qquad\qquad
    + [ A/ (\hat V_T^+ + gf\hat V_T^+ )] - [ \hat V_T / gf\hat V_T] + [\hat V_T / g\hat V_T]
    \end{aligned}
    $$}
Finally, apply the relations given by sequence $(7)$ to the second
term of the last expression and by sequence $(6)$ to the fourth
one
    {\small
    $$
    \begin{aligned}
    g_*[f] \,& =\,
    - [\hat V_T^++g\hat V_T^+ / g\hat V_T^+] - [ A / \hat V_T^+ ]
    + [ \hat V_T^++g\hat V_T^+ / \hat V_T^+ ]
    + [\hat V_T^++gf\hat V_T^+ / gf\hat V_T^+] +
     \\
    &\qquad\qquad
    + [ A / \hat V_T^+  ] - [\hat V_T^++gf\hat V_T^+ / \hat V_T^+]
    - [ \hat V_T / gf\hat V_T] + [\hat V_T / g\hat V_T]
    \,=
    \\
    & =\,
    - [\hat V_T^++g\hat V_T^+ / g\hat V_T^+]  + [ \hat V_T^++g\hat V_T^+ / \hat V_T^+]
    + [\hat V_T^++gf\hat V_T^+ / gf\hat V_T^+] -
     \\
    &\qquad\qquad
    - [\hat V_T^++gf\hat V_T^+ / \hat V_T^+]
    - [ \hat V_T / gf\hat V_T] + [\hat V_T / g\hat V_T]
    \,=
    \\
    & =\,
    [gf] - [g]
    \end{aligned}
    $$}
and the claim is proven.
\end{pf}

\begin{cor}{\,\hfill}
\begin{itemize}
    \item $[1]=0 \in K(S)$.
    \item $g_*[1]=[1]$.
    \item $1_*[g]=[g]$.
    \item $h_* g_* [f]=(hg)_*[f]$.
\end{itemize}
\end{cor}

Further,  $g$ yields a quasi-isomorphism
$M_{\bullet}\overset{g}\to g_* M_{\bullet}$ and, hence, an
identification $[M_{\bullet}]= [g_* M_{\bullet}]$ in $K(S)$, which
will be denoted by $g$. In particular, conjugation by it yields a
group isomorphism
    $$
    g^*\,:\, \aut(\Det([M_{\bullet}]))\,\overset{\sim}\longrightarrow\, \aut(\Det([g_* M_{\bullet}]))
    $$
which is closely related to the {\sl natural transformations} appearing in \cite{AP},~\S2.4.4.

\begin{defn}
Let us introduce
    $$
    \widehat{\mathcal Q}_S \,:=\,
    \left\{
    \begin{gathered}
        \text{pairs $(f,\tilde f)$ such that $f\in \bar{\mathcal Q}_S$ and}
        \\
        \text{$\tilde f$ is an automorphism of }\Det([f])
    \end{gathered}
    \right\}
    $$
\end{defn}

Given $(g,\hat g) , (f,\hat f)$, two elements of $\widehat{\mathcal
Q}_S$, observe that $g^*(\hat f)\in \aut(\Det(g_*[f]))$ since
$\hat f\in \aut(\Det([f]))$. Therefore, one has that $g^*(\hat
f)\otimes \hat g$ belongs to $\aut(\Det(g_*[f])\otimes\Det(
[g]))\simeq \aut(\Det([gf]))$. Summing up, one has

\begin{prop}
The composition law
    $$
    (g,\hat g) * (f,\hat f) \,:=\,
    (gf, g^*(\hat f)\otimes \hat g)
    $$
endows $\widehat{\mathcal Q}_S$ with a monoid structure such that the forgetful map
    $$
    \widehat{\mathcal Q}_S \,\longrightarrow\, \bar{\mathcal Q}_S
    $$
is a morphism of monoids.
\end{prop}

Further, since all the previous constructions are compatible with restrictions to non-empty open subsets of $S$, they can be stated in terms of sheaves on $S$. Observing that the very definitions of ${\mathcal Q}_S$ and $\bar{\mathcal Q}_S$ imply that
    $${\mathcal Q}_S \,\hookrightarrow\, \bar{\mathcal Q}_S$$
the following theorem is straightforward.

\begin{thm}\label{thm:centralextensionbarQ}
There is a central extension of sheaves on monoids on $S$
    $$
    0 \to {\mathbb G}_{m,S} \longrightarrow
    \widehat{\mathcal Q}_S \longrightarrow \bar{\mathcal Q}_S \to 0
    $$
that, when restricted to ${\mathcal Q}_S$, coincides with the extension of Theorem~\ref{thm:factorsystemQsheaf}.
\end{thm}

Let us finish with a sketch of an arithmetic application of our results that will be implemented in a future paper. The idea is to apply this machinery to study analogues of reciprocity laws in the case of families of curves $C\to S$ (including the non-abelian case) defined over ${\mathbb Z}$ and over $p$-adic numbers.

Note that, as a corollary of the above theorem, the factor system
    $$
    \bar{\mathcal Q}_S \times \bar{\mathcal Q}_S \, \overset{\bar \alpha_S}\longrightarrow \, {\mathbb G}_{m,S}
    $$
associated with the extension $\widehat{\mathcal Q}_S$  (Theorem~\ref{thm:centralextensionbarQ}) is equivalent to the factor system of $\tilde{\mathcal Q}_S$ (Theorem~\ref{thm:factorsystemQsheaf}) when restricted to ${\mathcal Q}_S \times {\mathcal Q}_S$.

However, what is even more relevant is that for any $g\in \bar{\mathcal Q}_S(T)$ there exists an open subset $R\subseteq T\subseteq S$ such that $g\vert_R\in {\mathcal Q}_S(R)$ (see Remark~\ref{rem:conditonsThmOpen}). In particular, the stalks of ${\mathcal Q}_S $ and $ \bar{\mathcal Q}_S$ at the generic point of $S$ coincide.

Finally, explicit expressions for the cases of curves $C$ over a field and over local artinian rings are well known (\cite{C,AP,MPa}). These arguments will be the key points for obtaining explicit expressions for the factor system $\bar\alpha_S$.


\end{document}